# LIMIT THEOREMS FOR THE TYPICAL POISSON–VORONOI CELL AND THE CROFTON CELL WITH A LARGE INRADIUS

By Pierre Calka and Tomasz Schreiber[1]

*Université René Descartes Paris 5 and Nicolaus Copernicus University*

In this paper, we are interested in the behavior of the typical Poisson–Voronoi cell in the plane when the radius of the largest disk centered at the nucleus and contained in the cell goes to infinity. We prove a law of large numbers for its number of vertices and the area of the cell outside the disk. Moreover, for the latter, we establish a central limit theorem as well as moderate deviation type results. The proofs deeply rely on precise connections between Poisson–Voronoi tessellations, convex hulls of Poisson samples and germ–grain models in the unit ball. Besides, we derive analogous facts for the Crofton cell of a stationary Poisson line process in the plane.

**1. Introduction and main results.** Consider $\Phi = \{x_n; n \geq 1\}$ a homogeneous Poisson point process in $\mathbb{R}^2$, with the two-dimensional Lebesgue measure $V_2$ for intensity measure. The set of cells

$$C(x) = \{y \in \mathbb{R}^2; \|y - x\| \leq \|y - x'\|, x' \in \Phi\}, \qquad x \in \Phi$$

(which are almost surely bounded polygons) is the well-known *Poisson–Voronoi tessellation* of $\mathbb{R}^2$. Introduced by Meijering [17] and Gilbert [8] as a model of crystal aggregates, it provides now models for many natural phenomena such as image analysis [18], molecular biology [7], thermal conductivity [15] and telecommunications [1, 2]. An extensive list of the areas in which the tessellation has been used can be found in [32] and [22].

In order to describe the statistical properties of the tessellation, the notion of *typical cell* $\mathcal{C}$ in the Palm sense is commonly used [20]. Consider the space $\mathcal{K}$ of convex compact sets of $\mathbb{R}^2$ endowed with the usual Hausdorff metric.

Received February 2004; revised July 2004.

[1]Supported in part by the Foundation for Polish Science (FNP).

*AMS 2000 subject classifications.* Primary 60D05, 60F10; secondary 60G55.

*Key words and phrases.* Germ–grain models, extreme point, large and moderate deviations, Palm distribution, Poisson–Voronoi tessellation, random convex hulls, stochastic geometry, typical cell.

Let us fix an arbitrary Borel set $B \subset \mathbb{R}^2$ such that $0 < V_2(B) < +\infty$. The distribution of the typical cell $\mathcal{C}$ is determined by the identity [20]:

$$\mathbf{E}\, h(\mathcal{C}) = \frac{1}{V_2(B)} \mathbf{E} \sum_{x \in B \cap \Phi} h(C(x) - x),$$

where $h : \mathcal{K} \longrightarrow \mathbb{R}$ runs throughout the space of bounded measurable functions.

Consider now the cell

$$C(0) = \{y \in \mathbb{R}^2; \|y\| \le \|y - x\|, x \in \Phi\}$$

obtained when the origin is added to the point process $\Phi$. It is well known [20] that $C(0)$ and $\mathcal{C}$ are equal in law. From now on, we will use $C(0)$ as a realization of the typical cell $\mathcal{C}$.

Let us denote by $R_m$ (resp. $R_M$) the radius of the largest (resp. smallest) disk centered at the origin included in (resp. containing) $C(0)$ and by $D(x, r)$, $x \in \mathbb{R}^2$, $r > 0$, the closed disk centered at $x$ and of radius $r$. The boundary of the polygon $C(0)$ then is contained in the annulus $\mathcal{A} = D(0, R_M) \setminus D(0, R_m)$. In [4], an explicit formula for the joint distribution of the pair $(R_m, R_M)$ and a characterization of the asymptotic behavior of the tail of the law of $R_M$ given $R_m$ were obtained. In particular, it was proved that *conditioning on the event* $\{R_m = r\}$, $r > 0$, *the thickness of the annulus* $\mathcal{A}$ *is a.s. "of order* $r^{-1/3}$*" when* $r$ *goes to infinity* (Result A).

Besides, a recent work by Hug, Reitzner and Schneider [11] has provided a proof (valid for any dimension) of D. G. Kendall's conjecture: *the shape of the typical Poisson–Voronoi cell, given that the area of the cell goes to infinity, tends a.s. to a disk* (Result B). This last result is stronger than Result A in the sense that the conditioning only holds on the area and not on the inradius but it is also weaker because it does not give such precise estimates for the thickness of the smallest annulus containing the boundary of $\mathcal{C}$.

A natural question arising from Result A is: how to estimate precisely the growth of the number of vertices of $C(0)$ and the decrease of the area of $C(0)$ outside the indisk when the inradius goes to infinity?

Let us denote by $\mathcal{C}_r$ (resp. $N_r$) a random variable taking values in the space of compact convex sets of $\mathbb{R}^2$ endowed with the Hausdorff metric (resp. in $\mathbb{N}$) whose distribution is given by the law of $C(0)$ [resp. the number of sides of $C(0)$] conditioned by the event $\{R_m = r\}$. It is well known [19] that

$$(\Phi | R_m = r) \stackrel{D}{=} \Phi_r \cup \{(2r) \cdot X_0\},$$

where $\Phi_r$ is a Poisson point process of intensity measure $\mathbf{1}_{D(0,2r)^c}(x)\, dx$ and $X_0$ is a uniform point on the boundary of $\mathbb{D} = D(0, 1)$. The cell $\mathcal{C}_r$ is then equal in law up to a uniform random rotation to the zero cell (i.e., the cell containing the origin) of the line process consisting of the bisecting lines of the segments between zero and the points of the process $\Phi_r \cup \{2r \cdot x_0\}$, where $x_0$ is the deterministic point $(1;0)$. The processes $\Phi_r$ and related random objects can be coupled on a common probability space in several natural ways; the coupling which we shall have in mind without further mentioning whenever stating $L^1$ or a.s. convergence results below is constructed in (2).

In the present paper we focus our interest on the asymptotic behavior of $N_r$ and $V_2(\mathcal{C}_r \setminus D(0,r))$ when $r \to +\infty$. Explicit formulae for the distributions of the number of sides and the area of the typical Poisson–Voronoi cell have been recently obtained (see [5, 6]) but it seems difficult to use them to obtain asymptotic results.

Note that for a regular polygon of indisk $D(0,r)$ whose vertices are all located in $\partial D(0, r + r^{-1/3})$, the number of sides and the area outside the indisk are asymptotically equivalent to $(\pi/\sqrt{2})r^{2/3}$ and $(2\pi/3)r^{2/3}$, respectively. The intuition provided by the results of [4] on the thickness of the annulus $\mathcal{A} = D(0,R_M) \setminus D(0,R_m)$ conditioned on $\{R_m = r\}$, as discussed above, suggests that the conditioned cell $\mathcal{C}_r$ should have the number of its sides $N_r$ of the same asymptotic order $r^{2/3}$.

Our first main result states that this is indeed the case and the growth rate for $N_r$ is exactly of the anticipated order $r^{2/3}$.

THEOREM 1. *When $r \to +\infty$, we have:*

(i) $\mathbf{E}N_r \sim a_1 r^{2/3}$,
(ii) $\lim_{r\to+\infty} N_r (a_1 r^{2/3})^{-1} = 1$ *in* $L^1$,

*where* $a_1 = 4\pi \cdot 3^{-1/3}\Gamma(5/3) \approx 7.86565$.

Note that we write $\alpha_r \sim \beta_r$ to indicate that $\lim_{r\to\infty} \alpha_r/\beta_r = 1$.

REMARK 1. It should be emphasized that our proof of Theorem 1 relies on an asymptotic equivalence between $N_r$ and the number of vertices of the convex hull generated by a homogeneous Poisson point process of intensity $4r^2$, which we establish below.

This equivalence is easily verified to be strong enough to also conclude a central limit theorem and variance asymptotics of order $r^{2/3}$ for $N_r$, should the corresponding results hold for the convex hulls. Such results are stated in Groeneboom's work [10] and were also stated in the previous version of our paper, yet upon its submission we have learned from several independent sources that some of the proofs of [10] may contain errors (although we do not know of what nature). However, we explain more precisely in a remark following the proof of Theorem 1 how the central limit theorem for $N_r$ can be deduced from Groeneboom's results.

The second theorem of this paper characterizes the asymptotic behavior of the area $V_2(\mathcal{C}_r \setminus D(0,r))$ which is proved to be of order $r^{2/3}$ up to a multiplicative constant. The obtained results include a central limit theorem and a moderate deviation principle.

THEOREM 2. *With $r \to +\infty$, we have:*

(A1) $\lim_{r\to+\infty} \frac{V_2(\mathcal{C}_r\setminus D(0,r))}{2\pi(4\pi)^{-2/3} b_1 r^{2/3}} = 1$ *in $L^1$ and a.s. for* $b_1 := \Gamma(\frac{2}{3})(\frac{\pi}{2})^{2/3}3^{-1/3}$,

(A2) $\mathbf{Var}\, V_2(\mathcal{C}_r \setminus D(0,r)) \sim b_2 r^{2/3}$ *for some constant* $b_2 > 0$,

(A3) $\frac{V_2(\mathcal{C}_r\setminus D(0,r)) \,-\, \mathbf{E}V_2(\mathcal{C}_r\setminus D(0,r))}{\sqrt{\mathbf{Var}\, V_2(\mathcal{C}_r\setminus D(0,r))}} \xrightarrow{D} \mathcal{N}(0,1)$,

(A4) *for each $\eta > 0$ we have*

$$I(\eta) := -\limsup_{r\to\infty} r^{-2/3} \log \mathbf{P}(V_2(\mathcal{C}_r \setminus D(0,r)) \geq (1+\eta)\mathbf{E}V_2(\mathcal{C}_r \setminus D(0,r))) > 0$$

*and, moreover,* $\lim_{\eta\to\infty} I(\eta)/\eta = (4\pi)^{1/3} \cdot b_1$ *with $b_1$ as in* (A1),

(A5) *for arbitrarily large $L > 0$ and arbitrarily small $\eta > 0$ there exists $Q := Q(L,\eta)$ such that*

$$\mathbf{P}(V_2(\mathcal{C}_r \setminus D(0, r + Qr^{-1/3})) \geq \eta r^{2/3}) = O(\exp(-Lr^{2/3})),$$

(A6) *for each $\eta > 0$ we have*

$$\liminf_{r\to\infty} \frac{1}{\log r} \log(-\log \mathbf{P}(V_2(\mathcal{C}_r \setminus D(0,r)) \geq (1-\eta)\mathbf{E}V_2(\mathcal{C}_r \setminus D(0,r)))) \geq \frac{2}{3}.$$

We can likewise obtain limit theorems for the Crofton cell of a stationary Poisson line process (see in particular [9, 12, 13] about the limit shape of the Crofton cell with a large area). More precisely, let us consider $\Phi'$ a Poisson point process in $\mathbb{R}^2$ of intensity measure (in polar coordinates) $\mathbf{1}_{\mathbb{R}_+}(r)\, dr\, d\theta$. The line process associated with $\Phi'$ (which is invariant in law by any translation of the plane) consists of the set of lines

$$L(x) = \{y \in \mathbb{R}^2; \langle y - x, x\rangle = 0\}, \qquad x \in \Phi', \tag{1}$$

where $\langle\cdot,\cdot\rangle$ denotes the usual scalar product on $\mathbb{R}^2$. Denoting by $H(x)$ the half-plane delimited by $L(x)$, $x \in \mathbb{R}^2$, and containing the origin, the Crofton cell $\mathcal{P}_0$ is given by the equality

$$\mathcal{P}_0 = \bigcap_{x\in\Phi'_t} H(x).$$

We successively define the radius $R'_m$ of the largest disk centered at the origin included in $\mathcal{P}_0$, $\mathcal{P}_r$ a random polygon distributed as the Crofton cell $\mathcal{P}_0$ conditioned by $\{R'_m = r\}$ and $N'_r$ the number of vertices of $\mathcal{P}_r$. In [4], we proved (see [4], Theorem 10) that when $r \to +\infty$, the boundary of $\mathcal{P}_r$ is included with "a great probability" in an annulus centered at the origin of thickness $r^{1/3}$. If the polygon $\mathcal{P}_r$ were regular, the number $N'_r$ and the area outside the indisk would be respectively of order $(\pi/\sqrt{2})r^{1/3}$ and $(2\pi/3)r^{4/3}$. Up to multiplicative constants, we are going to show that these orders are correct for both the number of vertices and the area outside the indisk.

The following two theorems are the equivalents of Theorems 1 and 2 for the Crofton cell $\mathcal{P}_0$.

THEOREM 3. *When $r \to +\infty$, we have:*

(i) $\mathbf{E}N'_r \sim a'_1 r^{1/3}$,
(ii) $\lim_{r\to+\infty} N'_r(a'_1 r^{1/3})^{-1} = 1$ *in* $L^1$,

*where* $a'_1 = 2^{4/3}\pi \cdot 3^{-1/3}\Gamma(5/3) \approx 4.95505$.

THEOREM 4. *With $r \to +\infty$, we have:*

(A1′) $\lim_{r\to+\infty} \frac{V_2(\mathcal{P}_r \setminus D(0,r))}{2\pi\pi^{-2/3}b_1 r^{4/3}} = 1$ *in* $L^1$ *and a.s. with* $b_1$ *as in Theorem* 2,
(A2′) $\mathbf{Var}\, V_2(\mathcal{P}_r \setminus D(0,r)) \sim b'_2 r^{7/3}$ *for some constant* $b'_2 > 0$,
(A3′) $\frac{V_2(\mathcal{P}_r \setminus D(0,r)) - \mathbf{E}V_2(\mathcal{P}_r \setminus D(0,r))}{\sqrt{\mathbf{Var}\, V_2(\mathcal{P}_r \setminus D(0,r))}} \xrightarrow{D} \mathcal{N}(0,1)$,
(A4′) *for each* $\eta > 0$ *we have*

$$\widetilde{I}(\eta) := -\limsup_{r\to\infty} r^{-1/3} \log \mathbf{P}(V_2(\mathcal{P}_r \setminus D(0,r)) \geq (1+\eta)\mathbf{E}V_2(\mathcal{P}_r \setminus D(0,r))) > 0$$

*and, moreover,* $\lim_{\eta\to\infty} \widetilde{I}(\eta)/\eta = \pi^{1/3} b_1$,

(A5′) *for arbitrarily large* $L > 0$ *and arbitrarily small* $\eta > 0$ *there exists* $Q := Q(L,\eta)$ *such that*

$$\mathbf{P}(V_2(\mathcal{P}_r \setminus D(0, r + Qr^{1/3})) \geq \eta r^{4/3}) = O(\exp(-Lr^{1/3})),$$

(A6′) *for each* $\eta > 0$ *we have*

$$\liminf_{r\to\infty} \frac{1}{\log r} \log(-\log \mathbf{P}(V_2(\mathcal{P}_r \setminus D(0,r)) \geq (1-\eta)\mathbf{E}V_2(\mathcal{P}_r \setminus D(0,r)))) \geq 1/3.$$

The methods for proving Theorems 1 and 2 on the one hand and Theorems 3 and 4 on the other hand are very similar so from now on we will essentially concentrate on the Poisson–Voronoi typical cell $\mathcal{C}$.

Instead of taking a limit when the value $r$ of the inradius goes to infinity, we shall rewrite the number $N_r$ so that the asymptotic results will be obtained as the intensity of the underlying Poisson point process in the plane goes to infinity. The area $V_2(\mathcal{C}_r \setminus D(0,r))$ will be dealt with along the same lines.

To this end, let us denote by $\Lambda$ the Poisson point process on $\mathbb{R}^2 \times \mathbb{R}_+$ of intensity measure $\mathbf{1}_{\mathbb{D}^c}(x)\mathbf{1}_{\mathbb{R}_+}(t)\, dx\, dt$ and by $\Psi_t$, $t \geq 0$, the Poisson point process on $\mathbb{R}^2$ defined by

$$\Psi_t = \{x \in \mathbb{R}^2 \setminus \mathbb{D}; \exists\, s \leq t, (x,s) \in \Lambda\}. \tag{2}$$

Then $r^{-1}\mathcal{C}_r = h_{1/r}(\mathcal{C})$ [where $h_{1/r}(x) = (1/r)\cdot x$, $x \in \mathbb{R}^2$] is easily seen to coincide in law with the zero cell $C_0^r$ of the line process consisting of the set of lines $L(x)$, $x \in \Psi_{4r^2} \cup \{x_0\}$ [see (1) and recall that $x_0$ stands for the deterministic point $(1;0)$ as defined above]. In other words, we have

$$C_0^r = \bigcap_{x\in\Psi_{4r^2}} H(x) \cap H(x_0). \tag{3}$$

In particular, $N_r$ coincides with the number of vertices of $C_0^r$ while $V_2(\mathcal{C}_r \setminus D(0,r))$ has the same law as $r^2 V_2(C_0^r \setminus \mathbb{D})$. Therefore, throughout the paper we will study the asymptotic behavior of $C_0^r$ rather than directly that of $\mathcal{C}_r$. Note also that, as already mentioned above, the relation (2) provides the coupling of the random objects considered in this paper on a common probability space, which we have in mind whenever stating a.s. or $L^1$-type results.

As already mentioned above, the proofs of Theorems 1 and 2 strongly rely on a connection, established via inversion of the complex plane, between the problems of determining the asymptotics of $N_r$ and $V_2(\mathcal{C}_r \setminus D(0,r))$, and some results on the asymptotic behavior of convex hulls of high-intensity Poisson point processes inside the disk $\mathbb{D}$, existing in the literature (see [3, 16, 24, 25, 29, 30]).

An extension of our results to higher dimensions will be given in a future paper.

## 2. Proofs.

In Lemma 1, we first relate $N_r$ to the number of vertices of the convex hull of a certain Poisson point process in $\mathbb{D}$ denoted by $Y_t$, $t = 4r^2$. Moreover, in Lemma 2 we represent the area $V_2(C_0^r \setminus \mathbb{D})$ as a defect measure of a certain germ–grain model in $\mathbb{D}$, generated by $Y_t$. Then, Lemmas 3 and 4 provide us with a comparison method between $Y_t$ and some homogeneous Poisson point processes so that the classical results on convex hulls due to Rényi and Sulanke [24] and Massé [16] can be applied, yielding a description of the asymptotic behavior of $N_r$ as stated in Theorem 1. The assertions of Theorem 2 are then concluded by combining the comparison Lemma 3 with appropriate results in [3, 28, 29, 30, 31].

LEMMA 1. *For every $r > 0$, $N_r$ coincides in law with the number $\widehat{N}_{4r^2}$ of vertices of the convex hull generated by the process $Y_{4r^2} \cup \{x_0\}$, where $Y_t$, $t \ge 0$, is a Poisson point process inside the disk $\mathbb{D}$, of intensity measure (in polar coordinates) $t \cdot \mu(d\rho, d\theta) = (t/\rho^3)\mathbf{1}_{(0,1)}(\rho)\, d\rho\, d\theta$.*

PROOF. Let us consider the inversion $I$ on $\mathbb{R}^2 \setminus \{0\}$ defined by

$$I(x) = \frac{1}{\|x^2\|} \cdot x, \qquad x \neq 0.$$

$I$ is a well-known involutive application which preserves the boundary of $\mathbb{D}$, transforms the interior of $\mathbb{D}$ into its exterior and conversely. In particular, $I$ transforms any line or circle into a line or a circle.

An easy calculation shows that the image in $I$ of the process $\Psi_t$, $t \geq 0$, is an inhomogeneous but rotation-invariant Poisson point process $Y_t$ in $\mathbb{D}$, of intensity measure $(t/r^3)\mathbf{1}_{(0,1)}(r)\,dr\,d\theta$ to be denoted by $t \cdot \mu$. Any line $L(x)$, $x \in \Psi_t$, is transformed into a circle having the segment $[0; I(x)]$ as its diameter. Note also that $I(x_0) = x_0$. Let us denote by $G(y) := D(y/2, \|y\|/2)$ a position-dependent grain and by

$$Y^{[t]} = \bigcup_{y \in Y_t} G(y)$$

the germ–grain model in $\mathbb{D}$ associated with the process $Y_t$ (see [30]). The image of the complement of the disk $\mathbb{D}$ in the Crofton cell $C_0^r$ is then the set $\mathbb{D} \setminus [Y^{[t]} \cup G(x_0)]$. Consequently, we have that

$$N_r = \#\{y \in Y_{4r^2} \cup \{x_0\}; \partial G(y) \cap \partial [Y^{[4r^2]} \cup G(x_0)] \neq \varnothing\}. \tag{4}$$

Let us notice that the boundary of the grain $G(z)$ associated with a given point $z \in Y_t$ intersects the boundary of the union of grains $Y^{[t]} \cup G(x_0)$ if and only if the convex hulls of $Y_t \cup \{x_0\}$ and of $(Y_t \setminus \{z\}) \cup \{x_0\}$ have different support functions. This yields the equivalence

$$\partial G(z) \cap \partial \left[ \bigcup_{y \in Y_t \cup \{x_0\}} G(y) \right] \neq \varnothing$$

$$\iff \quad (z \text{ is an extremal point of the convex hull of } Y_t \cup \{x_0\}).$$

Equality (4) implies then that $N_r$ is precisely the number $\widehat{N}_{4r^2}$ of points on the boundary of the convex hull of the process $Y_{4r^2} \cup \{x_0\}$. $\square$

The lemma below is a direct conclusion of the proof of Lemma 1.

Lemma 2. *With the notation as in the proof of Lemma 1, for each $r > 0$ the area $V_2(C_0^r \setminus \mathbb{D})$ coincides in distribution with the measure $\mu(\mathbb{D} \setminus [Y^{[t]} \cup G(x_0)])$ for $t = 4r^2$.*

Let us consider for every $\alpha > 0$ and $t \geq 1$ the event

$$A_{t,\alpha} = \{D(0, 1 - 2^{3\alpha} t^{-\alpha}) \not\subset Y^{[t]} \cup G(x_0)\}. \tag{5}$$

The following lemma shows that both the vertex process of the convex hull of $Y_t$ and the defect measure $\mu(\mathbb{D} \setminus Y^{[t]})$ are concentrated with an overwhelming probability in a close vicinity of the boundary $\partial\mathbb{D}$.

LEMMA 3. *There exists a positive constant $c>0$ such that for every $0\le\alpha<2/3$:*

(i) $\mathbf{P}(A_{t,\alpha}) \underset{t\to+\infty}{=} O(e^{-c\cdot t^{(1-(3/2)\alpha)}})$ *and* $\lim_{t\to\infty}\mathbf{1}_{A_{t,\alpha}}=0$ *a.s.;*

(ii) $\mathbf{E}(\widehat{N}_t\mathbf{1}_{A_{t,\alpha}}) \underset{t\to+\infty}{=} o(t^{1/3})$*;*

(iii) $t\mathbf{E}(\mu(\mathbb{D}\setminus Y^{[t]})\mathbf{1}_{A_{t,\alpha}}) \underset{t\to+\infty}{=} o(t^{1/3})$.

PROOF. (i) Applying the inversion $I$, we get from the equality $I(Y^{[t]}\cup G(x_0))=[C_0^{\sqrt{t}/2}]^{\mathbf{c}}$ that for every $\alpha\ge 0$,

$$\{D(0,1-2^{3\alpha}t^{-\alpha})\not\subset Y^{[t]}\cup G(x_0)\}=\{D(0,(1-2^{3\alpha}t^{-\alpha})^{-1})\not\supset C_0^{\sqrt{t}/2}\}. \tag{6}$$

The asymptotic result ([4], Theorem 5) on the distribution of the radius $R_M$ conditioned by the value $r$ of the inradius $R_m$ can be rewritten as: for every $0<c<8/(3\sqrt{2})$ and $0<\alpha'<1/3$,

$$\mathbf{P}\{D(0,r+r^{-\alpha'})\not\supset\mathcal{C}_r\} \underset{r\to+\infty}{=} O(e^{-c\cdot r^{(1-3\alpha')/2}}). \tag{7}$$

Since $C_0^r$ is the scale-$1/r$ homothetic image of the cell $\mathcal{C}_r$, we deduce from (7) that

$$\mathbf{P}\{D(0,1+r^{-(\alpha'+1)})\not\supset C_0^r\} \underset{r\to+\infty}{=} O(e^{-c\cdot r^{(1-3\alpha')/2}}). \tag{8}$$

Replacing $r$ by $\sqrt{t}/2$ in the preceding result and combining the equality of events (6) with the inequality $1+2^{2\alpha}t^{-\alpha}\le(1-2^{3\alpha}t^{-\alpha})^{-1}$ for $t$ large enough and $\alpha>0$, we get the first assertion of (i). To obtain the almost sure convergence put

$$\tilde{A}_{t,\alpha}:=\{D(0,1-2^{3\alpha}\lceil t\rceil^{-\alpha})\not\subset Y^{[t]}\cup G(x_0)\}$$

with $\lceil\cdot\rceil$ standing for the (upper) integer value. Note that for $k\in\mathbb{N}$, $t\in(k-1,k]$, we have

$$\tilde{A}_{t,\alpha}\subseteq\{D(0,1-2^{3\alpha}k^{-\alpha})\not\subset Y^{[k-1]}\cup G(x_0)\},$$

with the probabilities of the right-hand side events easily verified to satisfy the same bound as that for $A_{t,\alpha}$ in the first part of (i). Consequently, applying the Borel–Cantelli lemma, we conclude that $\lim_{t\to\infty}\mathbf{1}_{\tilde{A}_{t,\alpha}}=0$. The proof of (i) is now completed by the observation that $A_{t,\alpha}\subseteq\tilde{A}_{t,\alpha}$ for all $t>0$.

(ii) Use the Hölder–Schwarz inequality to get

$$\mathbf{E}(\widehat{N}_t\mathbf{1}_{A_{t,\alpha}})\le\sqrt{\mathbf{P}(A_{t,\alpha})}\sqrt{\mathbf{E}\widehat{N}_t^2}. \tag{9}$$

Thus, in view of the assertion (i) it remains to show that

$$\mathbf{E}\widehat{N}_t^2=O(t^2). \tag{10}$$

By Lemma 1, the assertion (10) is equivalent to

$$\mathbf{E}{N_r}^2 = O(r^4),$$

which is proved by elementary arguments in the Appendix.

(iii) The proof goes much along the same lines as that of (ii) above. We use the Hölder–Schwarz inequality to get

$$\mathbf{E}(\mu(\mathbb{D} \setminus Y^{[t]})\mathbf{1}_{A_{t,\alpha}}) \leq \sqrt{\mathbf{E}\mu^2(\mathbb{D} \setminus Y^{[t]})}\sqrt{\mathbf{P}(A_{t,\alpha})}. \tag{11}$$

We shall show that

$$\mathbf{E}\mu^2(\mathbb{D} \setminus Y^{[t]}) = O(t^2) \tag{12}$$

which, in view of (11) and the assertion (i), is more than enough to establish (iii). Using Lemma 2, (12) is equivalent to

$$\mathbf{E}[V_2(\mathcal{C}_r \setminus D(0,r))]^2 = O(r^4).$$

The proof of this last result is postponed to the Appendix. □

Remark 2. It is clear that removing the extra deterministic grain $G(x_0)$ does not affect the validity of the above results. Indeed, recalling (3) and using a similar argument as for [4], Theorem 5, we obtain a result analogous to (8): for every $0 < \alpha < 1/3$,

$$\mathbf{P}\left\{D(0, 1 + r^{-(\alpha'+1)}) \not\supset \bigcap_{x \in \Psi_{4r^2}} H(x)\right\} \underset{r \to +\infty}{=} O(e^{-c \cdot r^{(1-3\alpha')/2}}).$$

It remains to adapt the proof of Lemma 3 in order to get from the preceding result that for every $0 \leq \alpha < 2/3$,

$$\mathbf{P}\{D(0, 1 - 2^{3\alpha} t^{-\alpha}) \not\subset Y^{[t]}\} \underset{t \to +\infty}{=} O(e^{-c \cdot t^{(2-3\alpha)/2}}). \tag{13}$$

Near the boundary of $\mathbb{D}$, the intensity measure of the process $Y_t$ is "not far" from a multiple of the Lebesgue measure. Let us denote by $X_t$ a homogeneous Poisson point process in $\mathbb{D}$ of intensity measure $t\mathbf{1}_{(0,1)}(\rho)\rho\, d\rho\, d\theta$, $t \geq 0$.

In the next lemma, we prove by a coupling method (in the spirit of [30], Lemma 2) that the trace of $Y_t$ in any annulus $\mathbb{D} \setminus D(0, 1 - \varepsilon)$, $0 < \varepsilon < 1$, can be seen as a superset of the trace of $X_t$ and a subset of the trace of $X_{t/(1-\varepsilon)^4}$.

Lemma 4. *For every $\varepsilon > 0$, there exists a coupling of the point processes $X_t$, $Y_t$ and $X_{t/(1-\varepsilon)^4}$ such that almost surely,*

$$X_t \cap [\mathbb{D} \setminus D(0, 1 - \varepsilon)] \subseteq Y_t \cap [\mathbb{D} \setminus D(0, 1 - \varepsilon)] \subseteq X_{t/(1-\varepsilon)^4} \cap [\mathbb{D} \setminus D(0, 1 - \varepsilon)].$$

PROOF. Consider a Poisson point process $\Pi$ on $\mathbb{D}\times\mathbb{R}_+$ with intensity measure $\mathbf{1}_{\mathbb{D}}(y)\mathbf{1}_{\mathbb{R}_+}(t)\,dy\,dt$. It is then easily verified that $X_t$ coincides in distribution with the set of points $\{y\in\mathbb{D};\exists\, s\le t,(y,s)\in\Pi\}$ and $X_{t/(1-\varepsilon)^4}$ with $\{y\in\mathbb{D};\exists\, s\le t/(1-\varepsilon)^4,(y,s)\in\Pi\}$. Likewise, $Y_t$ coincides in law with $\{y\in\mathbb{D};\exists\, s\le t/\|y\|^4,(y,s)\in\Pi\}$. Since every $y\in\mathbb{D}\setminus D(0,1-\varepsilon)$ satisfies

$$t\le\frac{t}{\|y\|^4}\le\frac{t}{(1-\varepsilon)^4},$$

these representations of the point processes $X_t, Y_t$ and $X_{t/(1-\varepsilon)^4}$ are easily seen to provide the required coupling. □

PROOF OF THEOREM 1. By Lemma 3(i), (ii), we have that $\mathbf{E}(\widehat{N}_t\mathbf{1}_{A_{t,\alpha}})=o(t^{1/3})$ for every $0<\alpha<2/3$. Consequently, it suffices to study the asymptotic behavior of the number $\widehat{N}_t$, $t=4r^2$, of vertices of the convex hull of $Y_t\cup\{x_0\}$ outside the event $A_{t,\alpha}$. We have on the event $A^{\mathbf{c}}_{t,\alpha}$ that the vertices of the convex hull of $Y_t\cup\{x_0\}$ are located in the annulus $\mathbb{D}\setminus D(0,1-2^{3\alpha}t^{-\alpha})$.

Let us denote by $M_t$ (resp. $\widetilde{M}_t$) the number of vertices of the convex hull of $X_t$ (resp. $X_t\cup\{x_0\}$). Applying Lemma 4 to $\varepsilon=2^{3\alpha}t^{-\alpha}$, we obtain on the event $A^{\mathbf{c}}_{t,\alpha}$ that any vertex of the convex hull of $Y_t\cup\{x_0\}$ (resp. $X_{t/(1-2^{3\alpha}t^{-\alpha})^4}$) either is a vertex of the convex hull of $X_t\cup\{x_0\}$ (resp. $Y_t\cup\{x_0\}$) or is a point of $(Y_t\setminus X_t)\cap[\mathbb{D}\setminus D(0,1-2^{3\alpha}t^{-\alpha})]$ (resp. $(X_{t/(1-2^{3\alpha}t^{-\alpha})^4}\setminus Y_t)\cap[\mathbb{D}\setminus D(0,1-2^{3\alpha}t^{-\alpha})]$).

Denoting by $R_t$ (resp. $S_t$) the number of points in $(Y_t\setminus X_t)\cap[\mathbb{D}\setminus D(0,1-2^{3\alpha}t^{-\alpha})]$ (resp. in $(X_{t/(1-2^{3\alpha}t^{-\alpha})^4}\setminus Y_t)\cap[\mathbb{D}\setminus D(0,1-2^{3\alpha}t^{-\alpha})]$), we then deduce the following inequalities (on the event $A^{\mathbf{c}}_{t,\alpha}$):

$$N_t\le\widetilde{M}_t+R_t \tag{14}$$

and

$$N_t\ge\widetilde{M}_{t/(1-t^{-\alpha})^4}-S_t. \tag{15}$$

It now comes from the coupling construction of the point processes in the annulus $\mathbb{D}\setminus D(0,1-t^{-\alpha})$ that $R_t$ and $S_t$ are Poisson variables of respective means

$$\begin{aligned}\mathbf{E}(R_t)&=t\mu(\mathbb{D}\setminus D(0,1-2^{3\alpha}t^{-\alpha}))-tV_2(\mathbb{D}\setminus D(0,1-2^{3\alpha}t^{-\alpha}))\\&=t\pi\frac{(2^{3\alpha+1}t^{-\alpha}-2^{6\alpha}t^{-2\alpha})^2}{(1-2^{3\alpha}t^{-\alpha})^2}\end{aligned}$$

and

$$\begin{aligned}\mathbf{E}(S_t)&=\frac{t}{(1-2^{3\alpha}t^{-\alpha})^4}V_2(\mathbb{D}\setminus D(0,1-2^{3\alpha}t^{-\alpha}))-t\mu(\mathbb{D}\setminus D(0,1-2^{3\alpha}t^{-\alpha}))\\&=t\pi\frac{(2^{3\alpha+1}t^{-\alpha}-2^{6\alpha}t^{-2\alpha})^2}{(1-2^{3\alpha}t^{-\alpha})^4}.\end{aligned}$$

For $\alpha \in (1/2, 2/3)$, we get that

$$R_t \text{ and } S_t \underset{t\to+\infty}{\longrightarrow} 0 \qquad \text{in mean.} \tag{16}$$

Consequently, using (14) and (15), it only remains to obtain the law of large numbers for $\widetilde{M}_t$. To this end, let us now compare the two quantities $M_t$ and $\widetilde{M}_t$, $t \geq 0$. Any vertex distinct from $x_0$ of the convex hull of $X_t \cup \{x_0\}$ is obviously a vertex of the convex hull of $X_t$. Conversely, let us denote by $p_t$ (resp. $q_t$) the point of $X_t$ located in the upper (resp. lower) half-disk of $\mathbb{D}$ such that there is no point of the point process $X_t$ above (resp. under) the line through $p_t$ (resp. $q_t$) and $x_0$. If such a point does not exist, we take $p_t = -x_0$ (resp. $q_t = -x_0$). Then any vertex of the convex hull of $X_t$ is either a vertex of the convex hull of $X_t \cup \{x_0\}$ or is discarded when we add $\{x_0\}$ to the set of points, that is, is a vertex of the convex hull of the points of $X_t$ located in the corner corresponding to $x_0$ of the quadrilateral $\mathcal{Q}_t$ with vertices $x_0, p_t, q_t, 0$.

Let us denote by $V_t$ the number of "discarded" vertices. Then we have

$$M_t + 1 - V_t \leq \widetilde{M}_t \leq M_t + 1. \tag{17}$$

Conditionally to the positions of $p_t$ and $q_t$, the distribution of the points of $X_t$ inside the quadrilateral $\mathcal{Q}_t$ constituted by $x_0, p_t, q_t, 0$ is the law of a homogeneous Poisson point process of intensity measure $t\mathbf{1}_{\mathcal{Q}_t}(x)\,dx$. Consequently, after making an affine transformation, the number $V_t$ is the number of vertices in the left-lower corner of the convex hull of a homogeneous Poisson point process of intensity $t$ in a square.

Using (5.1) in [25] (or equivalently Section 3 in [21]) and Corollary 1 in [16], we deduce that

$$\begin{aligned} &\lim_{t\to+\infty} 3\mathbf{E}(V_t)(2\ln t)^{-1} = 1, \\ &\lim_{t\to+\infty} 3V_t(2\ln t)^{-1} = 1 \qquad \text{in probability.} \end{aligned} \tag{18}$$

It remains to apply (5.2) in [25] and Corollary 2 in [16] in order to get that for $c_1 = (3\pi/2)^{-1/3}\Gamma(5/3)$

$$\begin{aligned} &\lim_{t\to+\infty} \mathbf{E}M_t(2\pi^{4/3}c_1t^{1/3})^{-1} = 1, \\ &\lim_{t\to+\infty} M_t(2\pi^{4/3}c_1t^{1/3})^{-1} = 1 \qquad \text{in probability.} \end{aligned} \tag{19}$$

Combining (17) with (18) and (19), we deduce an $L^1$-law of large numbers (a consequence of the convergence in probability combined with the convergence of the means) for $\widetilde{M}_t$ when $t \to +\infty$. Putting these conclusions together with the inequalities (14) and (15) and the convergence stated in (16), we obtain the required results of Theorem 1 for $N_t$, $t = 4r^2$ (with $a_1 = 2^{5/3}\pi^{4/3}c_1$). Note that even though some of the cited results were originally established for the binomial rather than Poisson samples, they admit straightforward modifications for the Poisson case as well, due to the fact that the asymptotic properties of $M_t$ as $t \to \infty$ are only affected by the behavior of the underlying sample in infinitesimally close neighborhoods of the boundary $\partial D$; see, for example, the Poisson approximation argument in Section 3, Lemma 3.2 of [10]. The proof is complete. $\square$

REMARK 3. In this remark we discuss a method of obtaining asymptotic variance estimates and the central limit theorem for the number of vertices $N_r$, provided Groeneboom's paper [10] is correct.

Using [10], relations (1.1) and (1.2), we get, with the same notation as in the proof of Theorem 1, that

$$\lim_{t\to+\infty} 27\,\mathbf{Var}(V_t)(10\ln t)^{-1} = 1. \tag{20}$$

Besides, applying [10], equality (1.3), Theorem 3.4, we obtain that there exists a positive constant $c_2$ such that

$$\lim_{t\to+\infty} \mathbf{Var}(M_t)(2\pi^{4/3}c_2)^{-1} = 1 \tag{21}$$

and

$$\frac{M_t - 2\pi^{4/3}c_1 t^{1/3}}{\sqrt{2\pi^{4/3}c_2 t^{1/3}}} \xrightarrow{D} \mathcal{N}(0,1). \tag{22}$$

Combining (17) with (20), (21) and (22), we deduce a central limit theorem for $\widetilde{M}_t$ when $t \to +\infty$. As for the law of large numbers, it remains to use the inequalities (14) and (15) and the convergence (16) to have that

$$\mathbf{Var}\, N_r \sim a_2^2 r^{2/3}$$

and

$$\frac{N_r - \mathbf{E}N_r}{\sqrt{\mathbf{Var}\, N_r}} \xrightarrow{D} \mathcal{N}(0,1),$$

where $a_2 = \sqrt{2^{5/3}\pi^{4/3}c_2}$.

REMARK 4. Reitzner has recently proved an almost sure convergence for the number of vertices of the convex hull of $X_t$ when $t \to +\infty$ [23]. However, his result is valid for unit-balls of dimension $d \geq 4$ so it cannot be applied in our context to obtain the almost sure convergence when $r \to +\infty$ for the number $N_r$. This last property requires some additional work on extreme points of homogeneous Poisson point processes that will take place in a future paper.

PROOF OF THEOREM 2. The proof uses the representation of $V_2(C_0^r \setminus \mathbb{D})$, and hence of $V_2(\mathcal{C}_r \setminus D(0,r)) \stackrel{d}{=} r^2 V_2(C_0^r \setminus \mathbb{D})$ in terms of the defect measure of a high-density germ–grain model in $\mathbb{D}$, as stated in Lemma 2. The area $V_2(\mathcal{C}_r \setminus D(0,r))$ coincides in distribution with $r^2\mu(\mathbb{D} \setminus [Y^{[4r^2]} \cup G(x_0)])$. Since the assertions of Theorem 2 are to be concluded from general results for high-density germ–grain models as stated in [28, 29, 30], the deterministic grain $G(x_0)$ stands as a nuisance and the first step of our proof is aimed at getting rid of this grain. To this end, we denote by $\rho_t$ the (random) radius of the largest disk $D(0,\rho_t)$ centered in 0, which is completely contained in $Y^{[t]}$, and we observe that, by standard geometry,

$$\mu(\mathbb{D} \setminus Y^{[t]}) - \mu(\mathbb{D} \setminus [Y^{[t]} \cup G(x_0)]) = O((1-\rho_t)^{3/2}). \tag{23}$$

Moreover, using the result (13) with $\alpha := 1/2$, we get

$$\mathbf{P}(\rho_t < 1 - \sqrt{t}^{-1}) = O(\exp(-ct^{1/4})), \qquad c > 0. \tag{24}$$

Putting (23) and (24) together we conclude that

$$\mathbf{P}(\mu(\mathbb{D} \setminus Y^{[t]}) - \mu(\mathbb{D} \setminus [Y^{[t]} \cup G(x_0)]) > t^{-3/4}) = O(\exp(-Ct^{1/4})) \tag{25}$$

for some positive constant $C$. Recalling that we set $t := 4r^2$, it is easily seen that (25) is more than enough to safely replace $\mu(\mathbb{D} \setminus [Y^{[t]} \cup G(x_0)])$ by $\mu(\mathbb{D} \setminus Y^{[t]})$ when proving the assertions of Theorem 2 below.

To proceed with our proof, we observe that our germ–grain model $Y^{[t]}$ in close neighborhoods of the boundary $\partial D$ "differs only negligibly" from the germ–grain model $X^{[\pi t]}$ as considered in Section 3 of [30], defined by $X^{[\pi t]} := \bigcup_{x \in X_t} G(x)$, where $X_t$ is the homogeneous Poisson point process of intensity $t$, restricted to $\mathbb{D}$ (see the notation introduced in the discussion preceding Lemma 4). Indeed, it follows by Lemma 4 that for arbitrarily small $\delta > 0$, taking $\varepsilon := 1 - (1+\delta)^{-1/4}$, we can find a coupling of versions of $X_t$, $X_{t(1+\delta)}$ and $Y_t$ such that almost surely

$$\bigcup_{y \in X_t \cap [\mathbb{D} \setminus D(0,1-\varepsilon)]} G(y) \subset \bigcup_{y \in Y_t \cap [\mathbb{D} \setminus D(0,1-\varepsilon)]} G(y) \subset \bigcup_{y \in X_{(1+\delta)t} \cap [\mathbb{D} \setminus D(0,1-\varepsilon)]} G(y).$$

Moreover, in view of Lemma 3(iii), applied with $\alpha := 1/2$ [see also (24)], for the purpose of the proof of Theorem 2 we can safely ignore the behavior of $Y^{[t]}$ inside $D(0, 1 - t^{-1/2})$.

A further observation in the same spirit is that in close neighborhoods of the boundary $\partial\mathbb{D}$, the measure $\mu$ here differs only negligibly from $2\pi\mu$ as considered in Section 3 of [30], to be denoted here by $2\pi\mu^*$ to avoid confusion and defined there by $2\pi\mu^*(dr, d\theta) = \mathbf{1}_{(0,1)}(r)\, dr\, d\theta$ in polar coordinates.

It is proved in Section 3 of [30] [see (18) there] that on the event that the convex hull of $X_t$ does contain the origin, the random variables $2\pi\mu^*(\mathbb{D} \setminus$ $X^{[\pi t]}$) and $\pi(2 - b_{\pi t})$ coincide, with $b_{\pi t}$ standing for the mean width of the convex hull generated by $X_t$. The probability of $\{0 \notin \mathrm{conv}(X_t)\}$ decays exponentially with $t$ in that there exists a constant $c > 0$ with $\mathbb{P}(0 \notin \mathrm{conv}(X_t)) = O(\exp(-ct))$, which is negligible in our setting; see, for example, Theorem 2 in [27] or (3.2) in [14].

This shows that when proving the assertions (A1)–(A6) of our theorem we can safely replace $\mu(\mathbb{D} \setminus Y^{[t]})$ by $\pi(2 - b_{\pi t}), t := 4r^2$ and, consequently, $V_2(\mathcal{C}_r \setminus D(0, r))$ by $\pi r^2(2 - b_{4\pi r^2})$ since $V_2(\mathcal{C}_r \setminus D(0, r))$ coincides in law with $r^2\mu(\mathbb{D} \setminus [Y^{[4\pi r^2]} \cup G(x_0)])$ and the effect of adding the extra deterministic grain $G(x_0)$ is negligible as discussed above. It puts us in a position to apply Theorem 2 in [30] combined with (20) in [30] (see also Theorem 6 in [28]), stating that $\lim_{t\to\infty} t^{2/3}\mathbb{E}(2 - b_{\pi t}) = 2b_1$ with $b_1$ as in (A1), to conclude that

$$\lim_{r\to\infty} \frac{\mathbb{E}V_2(\mathcal{C}_r \setminus D(0, r))}{2\pi(4\pi)^{-2/3}b_1 r^{2/3}} = 1.$$

The strong law of large numbers as stated in Corollary 2 and Corollary 3 in [29] allow us to conclude the assertion (A1) of our theorem (technically speaking, we get the convergence of means and the a.s. convergence, but these together yield immediately the required $L^1$ convergence).

The central limit theorem in (A3) follows now from Theorem 6 in [29]. Since for compact convex $K \subseteq \mathbb{R}^2$ we have the relation $b(K) = L(K)/\pi$, with $b$ standing for the mean width and $L$ for the perimeter (see page 210 in [26]), we could alternatively have used the results in Section 5 in [3], which yield also our assertion (A2) (see also Theorem 5 in [29]).

The assertions (A4) and (A5) are now direct consequences of Theorems 8 and 3 in [30]. The remaining assertion (A6) follows by Theorems 1 and 2 in [31].

Note that even though some of the cited results were originally established for the binomial rather than Poisson samples, they admit straightforward modifications for the Poisson case as well, due to the fact that the asymptotic properties of $\mu(\mathbb{D} \setminus Y^{[t]})$ as $t \to \infty$ are only affected by the behavior of the underlying germ point process in infinitesimally close neighborhoods of the boundary $\partial D$; see, for example, the comparison formulae (5), (6) and (9), (10) in [29]. The proof is complete. □

Proof of Theorems 3 and 4. The image of a Poisson point process of intensity measure $\mathbf{1}_{(r,+\infty)}(\rho)\,d\rho\,d\theta$, $r > 0$, by $I \circ h_{1/r}$ is a Poisson point process in the disk $\mathbb{D}$ of intensity measure (in polar coordinates) $r \cdot \nu(d\rho, d\theta) = (r/\rho^2)\mathbf{1}_{(0,1)}(\rho)\,d\rho\,d\theta$. Replacing the measure $\mu$ by $\nu$ in the preceding arguments, we easily obtain the results of Theorems 3 and 4. □

REMARK 5. In the same way as for the Poisson–Voronoi typical cell, we could prove a central limit theorem for the number $N'_r$, provided Groeneboom's results [10] are correct, that is,

$$\mathbf{Var}\, N'_r \sim a'^2_2 r^{1/3}$$

and

$$\frac{N'_r - \mathbf{E}N'_r}{\sqrt{\mathbf{Var}\, N'_r}} \xrightarrow{D} \mathcal{N}(0,1).$$

## APPENDIX

In this appendix, we give some technical results about the variables $N_r$ and $V_2(\mathcal{C}_r \setminus D(0,r))$ which are useful in the proof of Lemma 3. The proofs here only use elementary facts on the Poisson–Voronoi tessellation. In particular we do not need to use the analogy provided by Lemmas 1 and 2 with the convex hulls of the point processes inside the unit disk.

FACT 1. *There exist positive constants $K$ and $\lambda$ such that when $r \to +\infty$:*

(i) $\mathbf{E}(e^{\lambda V_2(\mathcal{C}_r \setminus D(0,r))}) = O(e^{Kr^2})$;
(ii) $\mathbf{E}(e^{\lambda N_r}) = O(e^{Kr^2})$;
(iii) *in particular, when $r \to +\infty$, we have* $\mathbf{E}({N_r}^2) = O(r^4)$ *and* $\mathbf{E}[(V_2(\mathcal{C}_r \setminus D(0,r)))^2] = O(r^4)$.

PROOF. (i) Let us apply the method provided by Gilbert [8] in order to estimate the expectation $\mathbf{E}(\exp\{V_2(\mathcal{C}_r \setminus D(0,r))\})$, $r > 0$.

We first define $R_{r,V}$ as the radius of the ball centered at the origin which has the same area as $\mathcal{C}_r \setminus D(0,r)$, $r > 0$. Then $R_{r,V}$ satisfies the following inequality for any $s \in (0,1)$:

$$\mathbf{E}\biggl(\int_{D(0,R_{r,V})} \exp(s\pi\|x\|^2)\,dx\biggr) \le \mathbf{E}\biggl(\int_{\mathcal{C}_r \setminus D(0,r)} \exp(s\pi\|x\|^2)\,dx\biggr). \tag{26}$$

Moreover,

$$\mathbf{E}\biggl(\int_{D(0,R_{r,V})} \exp(s\pi\|x\|^2)\,dx\biggr) = \frac{1}{s}\mathbf{E}(\exp\{sV_2(\mathcal{C}_r \setminus D(0,r))\}) - 1. \tag{27}$$

Recalling that $\mathcal{C}_r$ is up to a rotation equal in law to the zero cell delimited by the bisecting lines of the segments between the origin and the points of the process $\Phi_r \cup \{2r \cdot x_0\}$, we obtain

$$\mathbf{E}\biggl(\int_{\mathcal{C}_r \setminus D(0,r)} \exp(s\pi\|x\|^2)\,dx\biggr)$$

$$
\begin{aligned}
&= \int \mathbf{P}\{x \in \mathcal{C}_r \setminus D(0,r)\} \exp(s\pi\|x\|^2)\,dx \\
&\le \int \mathbf{P}\{\Phi_r \cap D(x,\|x\|) = \varnothing\} \exp(s\pi\|x\|^2)\,dx \qquad (28)\\
&\le 2\pi e^{4\pi r^2} \int_r^{+\infty} \exp((s-1)\pi r^2) r\,dr \\
&= e^{4\pi r^2} \cdot \frac{1}{1-s} e^{(s-1)\pi r^2}.
\end{aligned}
$$

Combining (26) with (27) and (28), we obtain the point (i) of Fact 1.

(ii) We apply the method due to Zuyev [33] to estimate the expectation $\mathbf{E}(\exp(sN_r))$, $r > 0$. Let $\mathcal{F}$ (resp. $\mathcal{I}$) be the union of the four open disks of radius 1 centered at the points $(\pm 1, 0)$, $(0, \pm 1)$ so that the origin lies on their boundary (resp. the set of points of $\mathcal{F}$ belonging to exactly two of these disks). Besides, we denote by $\mathcal{I}_1, \ldots, \mathcal{I}_4$ the connected components of $\mathcal{I}$.

As previously seen, $N_r$ coincides with the number of sides of the Crofton cell $C_0^r$ delimited by the lines $L(x)$ associated with the points $x \in \Psi_{4r^2} \cup \{x_0\}$ (see the Introduction). Suppose now that there exists $\alpha > 1$ such that $\Psi_{4r^2}$ intersects every connected component of $\alpha\mathcal{I}$. In that case, the number of edges of $C_0^r$ is at most equal to the number of points of $(\Psi_{4r^2} \cup \{x_0\}) \cap \alpha\mathcal{F}$. Consequently, writing $\mathcal{E}[\alpha]$ for the event that there exists at least one connected component of $\alpha\mathcal{I}$ which is not hit by $\Psi_{4r^2}$, we have that for any $\alpha > 1$ there exists $\delta, K > 0$ such that

$$
\begin{aligned}
&\mathbf{E}(\exp(sN_r - s)) \\
&\quad \le \mathbf{E}(\exp(s\#[\Psi_{4r^2} \cap \mathcal{F}])) \\
&\qquad + \sum_{n=1}^{+\infty} \mathbf{E}(\exp(s\#[\Psi_{4r^2} \cap \alpha^n \mathcal{F}]) \mathbf{1}_{\mathcal{E}[\alpha^{n-1}] \setminus \mathcal{E}[\alpha^n]}) \\
&\quad \le e^{4r^2(e^s-1)V_2(\mathcal{F})} \\
&\qquad + 4 \sum_{n=1}^{+\infty} \sum_{k=0}^{+\infty} e^{sk} \mathbf{P}\{\#[\Psi_{4r^2} \cap \alpha^n \mathcal{F} \cap (\alpha^{n-1}\mathcal{I}_1)^{\mathbf{c}}] = k\} \qquad (29)\\
&\qquad\qquad \times \mathbf{P}(\Psi_{4r^2} \cap \alpha^{n-1}\mathcal{I}_1 = \varnothing) \\
&\quad \le e^{4r^2(e^s-1)V_2(\mathcal{F})} \\
&\qquad + 4 \sum_{n=1}^{+\infty} \sum_{k=0}^{+\infty} \frac{e^{sk}}{k!} \{4r^2 V_2[(\alpha^n \mathcal{F}) \setminus (\alpha^{n-1}\mathcal{I}_1)]\}^k \cdot e^{-4r^2 V_2(\alpha^{n-1}\mathcal{I}_1 \setminus \mathbb{D})} \\
&\quad \le e^{4r^2(e^s-1)V_2(\mathcal{F})}
\end{aligned}
$$

$$+ 4e^{4\pi r^2} \sum_{n=1}^{+\infty} \exp\{4r^2\alpha^{2n-2}[(e^s - 1)(\alpha^2 v - w) - w]\},$$

where $v$ (resp. $w$) denotes the area of $\mathcal{F}$ (resp. $\mathcal{I}_1$).

When $[(e^s - 1)(\alpha^2 v - w) - w] < 0$ [i.e., $s < -\ln(1 - w/(\alpha^2 v))$] and $r > 1$, the series in (30) is convergent and bounded by a constant independent of $r$. Consequently, since $\alpha$ is arbitrarily chosen in $(1, +\infty)$, we have that for every $s < -\ln(1 - w/v)$ and $r > 1$, there exists $\delta, K > 0$ such that

$$\mathbf{E}(\exp(sN_r - s)) \le \delta e^{Kr^2}. \tag{30}$$

(iii) Applying Jensen's inequality to the convex function $a(x) = \exp(s\sqrt{x}\,)$, $s > 0$, $x \in [1/s^2, +\infty)$ and to the variable $\max(1/s^2, [V_2(\mathcal{C}_r \setminus D(0,r))]^2)$, we get

$$\exp\{s\sqrt{\mathbf{E}[(V_2(\mathcal{C}_r \setminus D(0,r)))^2]}\} \le \mathbf{E}(e^{s(V_2(\mathcal{C}_r \setminus D(0,r)))}) + e. \tag{31}$$

Combining (31) with the point (i), we deduce that when $r \to +\infty$

$$\mathbf{E}[(V_2(\mathcal{C}_r \setminus D(0,r)))^2] = O(r^4).$$

The same proof holds for $\mathbf{E}(N_r^2)$ as well. □

Remark 6. For the Crofton cell of a stationary Poisson line process, it is equally possible to use the same type of arguments to obtain that the second moments of the number of vertices and of the area of the complement of the indisk are at most of order $r^2$ when the inradius $r$ goes to infinity.

**Acknowledgment.** Special thanks are due to anonymous referees whose remarks have been very helpful in improving this paper.

Université René Descartes Paris 5
MAP5, UFR Math-Info
45, rue des Saints-Pères 75270
Paris Cedex 06
France
e-mail: pierre.calka@math-info.univ-paris5.fr

Nicolaus Copernicus University
ul. Chopina 12/18
87-100 Toruń
Poland
e-mail: tomeks@mat.uni.torun.pl